# ASYMPTOTIC BEHAVIOUR OF THE ONE-DIMENSIONAL "ROCK-PAPER-SCISSORS" CYCLIC CELLULAR AUTOMATON

By Benjamin Hellouin de Menibus

*Laboratoire de Recherche en Informatique, UMR 8623*
*Université Paris-Sud – CNRS – CentraleSupélec,*
*Université Paris-Saclay, France*

and

By Yvan Le Borgne
*Laboratoire Bordelais de Recherche en Informatique, UMR 5800*
*Université de Bordeaux – CNRS – Bordeaux INP,*
*F-33405 Talence, France*

The one-dimensional three-state cyclic cellular automaton is a simple spatial model with three states in a cyclic "rock-paper-scissors" prey-predator relationship. Starting from a random configuration, similar states gather in increasingly large clusters; asymptotically, any finite region is filled with a uniform state that is, after some time, driven out by its predator, each state taking its turn in dominating the region (heteroclinic cycles).

We consider the situation where each site in the initial configuration is chosen independently at random with a different probability for each state. We prove that the asymptotic probability that a state dominates a finite region corresponds to the initial probability of its prey. The proof methods are based on discrete probability tools, mainly particle systems and random walks.

Cyclic dominance is a general term for phenomena where different states (species, strategies, etc.) are in prey-predator relationships that form a cycle: A preys on B preys on C… preys on A. This phenomenon occurs in many natural or theoretical systems, among which a few examples are:

**Population ecology** male mating strategies in side-blotched lizard [30], antibiotic production and resistance in *E.Coli* [21], parasite-grass-forb interactions [6], oscillations in the population size of pacific salmon [15], etc.

**Game theory** pure or stochastic strategies in rock-paper-scissors type games, iterated prisoner's dilemma [29, 19], public goods games [16, 27], etc.

**Infection models** The SIRS compartmental model [1] (susceptible / infectious / recovered, when a recovered agent may become susceptible again), forest fire models [2], etc.

Many additional examples can be found in [32] (Section 7) and [33].

May and Leonard's [24] is the first effort to model the evolution of three species with cyclic dominance, using the standard Lotka-Volterra equations; it is a mean-field approximation, that is, it assumes the population is well-mixed. The system exhibits so-called *heteroclinic cycles* where each species in turn dominates almost the whole space before being replaced by its predator. Consequently, cyclic dominance has been proposed as a mechanism to explain the coexistence of various strategies or species [20] (biodiversity), the regular oscillations in population sizes of different species [15], and some counter-intuitive phenomena such as the "survival of the weakest" [14]. In other contexts, heteroclinic cycles appear to coincide with important concepts: for example, social choice among three cyclically dominant choices can lead to a heteroclinic cycle along the so-called *bipartisan set* [22].

Mean-field models do not take into account spatial aspects of the evolution of populations, such as the effect of population structure, mobility, dispersal, local survival, etc. This is why spatial models have been introduced both in ecology [7, 34] and in so-called evolutionary game theory [32]. In both cases agents have a spatial location and can only interact with their neighbours at short range. There is some variety in spatial models:

**Space** a lattice in one, two or more dimensions, or a graph with more structure;

**Updates** discrete or continuous time, synchronous or asynchronous updates;

**Dynamics** usually a predator replaces a prey by a copy of itself (*replicator dynamics*). The model can include empty space, different ranges, threshold effects, invasion probabilities, etc.;

**Boundaries** infinite, periodic or fixed boundary conditions, choice of the initial configuration.

In this article, we consider arguably the simplest spatial model for cyclic dominance: the one-dimensional, 3-state cyclic cellular automaton. To each site on the lattice $\mathbb{Z}$ is assigned an initial state in $\mathbb{Z}/3\mathbb{Z}$. At each (discrete) time step, every site is updated synchronously: if any of the two neighbouring sites contains a predator, it becomes the new state for this site.

While the restriction to one dimension may not be ecologically realistic (two-dimensional models being the object of more interest [34]), it has two benefits. First, its simple spatial structure makes many questions mathematically tractable, while the two-dimensional models have much more complex dynamics with structured interfaces between regions [9]. Second, its dynamics is similar to an *interacting particle system* with borders progressing at constant speed and annihilating on contact (*ballistic annihilation* - see Figure 1); this is a subject of independent interest [5] and many tools have been developed for it [3].

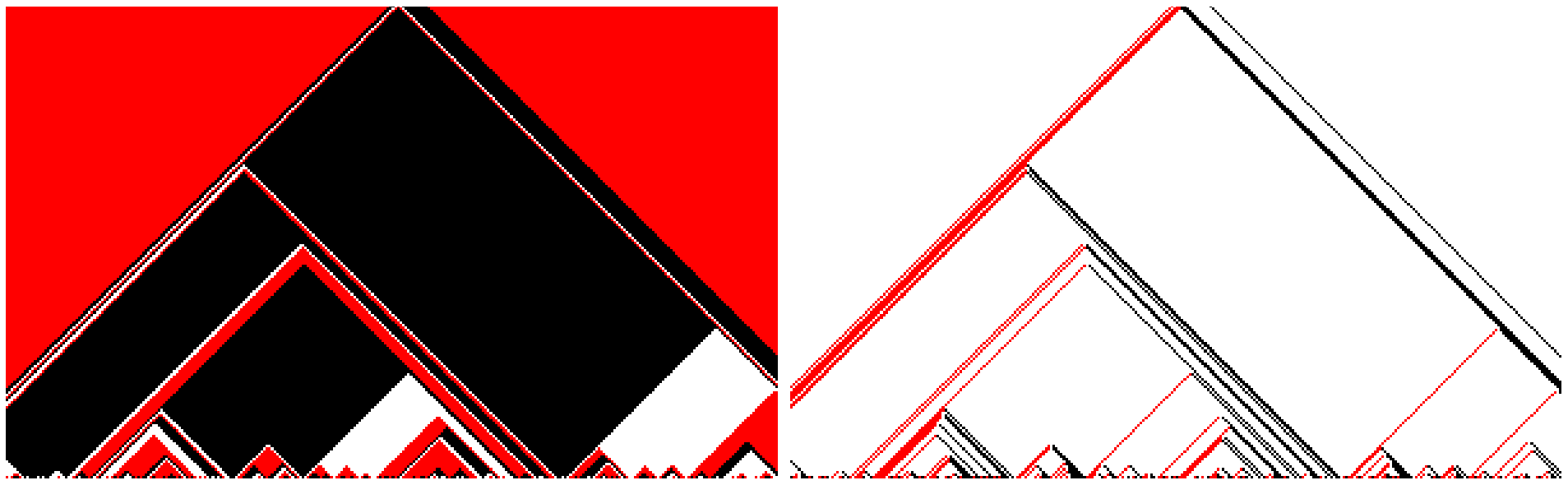

*Fig 1. (Left) The 3-state cyclic cellular automaton; (Right) The dynamics of its particles.*

*Note.* In all space-time diagrams of this article, the initial configuration is drawn horizontally at the bottom and time goes from bottom to top. States are represented by colours following the convention $0 \mapsto \square, 1 \mapsto \blacksquare, 2 \mapsto$ ■, $3 \mapsto$ ■, $4 \mapsto$ ■.

The seminal work of Fisch [10] focused on the case where each site is independently assigned a random state with uniform probability. He proved a clustering phenomenon: for 3 or 4 states, large monochromatic regions emerge and grow, but each region keeps changing state arbitrarily late (*fluctuation*, the spatial counterpart of heteroclinic cycle); for 5 states or more, the regions reach a limit size then stay unchanged (*fixation*). These behaviours, illustrated in Figure 2, are considered as a prime example of self-organisation in a relatively simple model [28]. These results were later refined in terms of cluster growth rate, number of state changes, etc. [11, 13, 23].

The present article focuses on the asymptotic behaviour of the 3-state cyclic cellular automaton when the initial configuration is chosen independently at random, but with distinct probabilities for each state, breaking the symmetry. It is not hard to see that the same clustering phenomenon as in the uniform case occurs. Our main result (Theorem 6) is that the asymptotic probability for any region to be dominated by a given state corresponds to the initial probability of its prey; this completely determines the limit probability measure. A similar relationship was observed empirically between invasion rates and asymptotic probability in more complex models [35]; see [32], Section 7.7 for a detailed account. However, we could not find a conjecture for this phenomenon in such a simple model, and this is the first formal proof of a similar result to our knowledge.

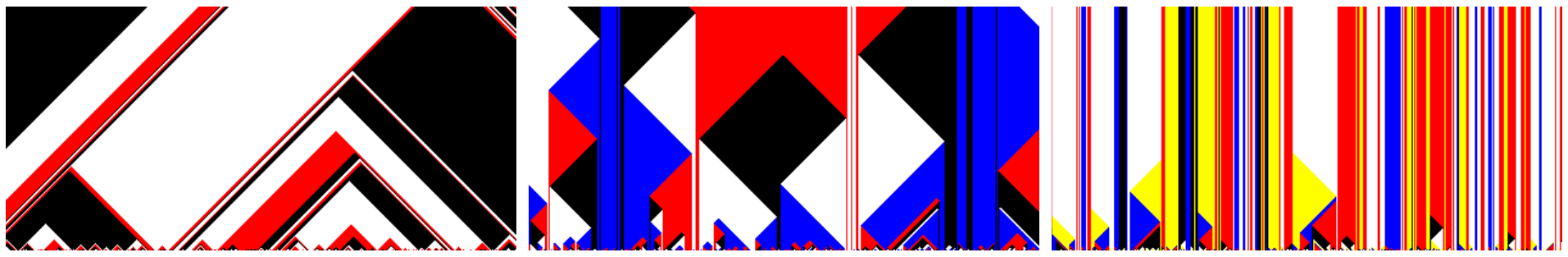

*Fig 2. Left to right: the 3-, 4- and 5-state cyclic cellular automata iterated on an initial configuration where each site is drawn in a uniform i.i.d. manner.*

Our approach is based on a correspondence between the time evolution of the borders and some well-chosen random walk, a method that was already used in the study of one-dimensional cellular automata [3]. Compared with previous work, the random walk is not the standard symmetric walk and the probability of a step up or down depends on the current position.

## 1. Definitions.

1.1. *Symbolic space.* For $\mathcal{A}$ a finite alphabet, define $\mathcal{A}^* = \bigcup_{n\in\mathbb{N}} \mathcal{A}^n$ the set of finite *patterns* (or words) and $\mathcal{A}^\mathbb{Z}$ the set of (one-dimensional) *configurations*, that is, the set of bi-infinite words over the alphabet $\mathcal{A}$. For example, for $a \in \mathcal{A}$, denote ${}^\infty a^\infty \in \mathcal{A}^\mathbb{Z}$ by ${}^\infty a^\infty{}_i = a$ for all $i \in \mathbb{Z}$. We endow $\mathcal{A}^\mathbb{Z}$ with the product topology of the discrete topology on $\mathcal{A}$.

For $u \in \mathcal{A}^*$, denote $|u|$ its length, and for $i \in \mathbb{Z}$, define the *cylinder* $[u]_i = \{x \in \mathcal{A}^\mathbb{Z} : x_{[i,i+|u|-1]} = u\}$, with $[u] = [u]_0$. Cylinders form a clopen basis of $\mathcal{A}^\mathbb{Z}$. A word $u \in \mathcal{A}^*$ is a *factor* of a configuration $x \in \mathcal{A}^\mathbb{Z}$ if $x \in [u]_i$ for some $i \in \mathbb{Z}$.

Define the *shift function* $\sigma : \mathcal{A}^\mathbb{Z} \to \mathcal{A}^\mathbb{Z}$ by $\sigma(x)_i = x_{i-1}$ for any $i \in \mathbb{Z}$. A *cellular automaton* is a pair $(\mathcal{A}, F)$ where $F : \mathcal{A}^\mathbb{Z} \to \mathcal{A}^\mathbb{Z}$ is a function that is continuous for the product topology and commutes with $\sigma$ (i.e. $F\circ\sigma = \sigma\circ F$). Alternatively, by the Curtis–Hedlund–Lyndon theorem [17], $F$ is defined by a finite *neighbourhood* $N \subset \mathbb{Z}$ and a local rule $f : \mathcal{A}^N \to \mathcal{A}$ in the sense that $F(x)_i = f((x_{i+j})_{j\in N})$.

In the figures, we represent the time evolution of cellular automata starting from an initial configuration $x \in \mathcal{A}^\mathbb{Z}$ by the corresponding two-dimensional *space-time diagram* $(F^t(x)_i)_{t\in\mathbb{N},\ i\in\mathbb{Z}}$.

The *frequency* of a finite word $u$ in a configuration $x \in \mathcal{A}^{\mathbb{Z}}$ is defined as:

$$\mathrm{freq}(u,x) = \limsup_{n\to\infty} \frac{1}{(2n+1)} \operatorname{Card}\{i \in \{-n,\ldots,n\} \ : \ x \in [u]_i\}.$$

1.2. *Cyclic cellular automata.*

Definition 1 ($n$-state cyclic cellular automaton). *$(\mathbb{Z}/n\mathbb{Z}, C_n)$ is the $n$-state cyclic cellular automaton defined on the neighbourhood $N = \{-1,0,1\}$ by the local rule $c_n$:*

$$c_n(u_{-1},u_0,u_1) = \begin{cases} u_0+1 & \text{if } u_1 = u_0+1 \bmod n \quad \text{or} \quad u_{-1} = u_0+1 \bmod n, \\ u_0 & \text{otherwise.} \end{cases}$$

All operations concerning $n$-state cyclic automata are assumed to be modulo $n$.

As should be clear from Figure 2, the self-organisation is driven by borders between monochromatic regions behaving as particles. We call particles the factors $ab$ of length 2 (with $a \neq b$) in a configuration. Each particle moves "from predator to prey", that is, left if $b = a+1$, right if $b = a-1$, and stays put otherwise. This motivates the following definitions:

**Positive particles** $p_+ = \{ab : b = a-1\}$;
**Negative particles** $p_- = \{ab : b = a+1\}$;
**Neutral particles** $p_= = \{ab : b \notin \{a-1,\ a,\ a+1\}\}$.

We write $[p_+]_i$ as a shorthand for $\bigcup_{ab\in p_+}[ab]_i$: it means that a positive particle occurs at position $i$. Define similarly $[p_=]_i$ and $[p_-]_i$. Notice that $p_= = \emptyset$ for $n = 3$. Figure 1 illustrates the particle dynamics for $n = 3$.

1.3. *Probability measures on $\mathcal{A}^{\mathbb{Z}}$.* Let $\mathfrak{B}$ be the Borel sigma-algebra of $\mathcal{A}^{\mathbb{Z}}$. Denote by $\mathcal{M}(\mathcal{A}^{\mathbb{Z}})$ the set of probability measures on $\mathcal{A}^{\mathbb{Z}}$ defined on the sigma-algebra $\mathfrak{B}$. Since the cylinders $\{[u]_n \ : \ u \in \mathcal{A}^*, n \in \mathbb{Z}\}$ form a basis of the product topology on $\mathcal{A}^{\mathbb{Z}}$, a measure $\mu \in \mathcal{M}(\mathcal{A}^{\mathbb{Z}})$ is entirely characterised by the values $\mu([u]_n)$.

In this paper, we only consider the set $\mathcal{M}_\sigma(\mathcal{A}^{\mathbb{Z}})$ of $\sigma$-invariant probability measures, and therefore write $\mu([u])$ instead of $\mu([u]_n)$.

*Examples.*

**Monochromatic measures.** For $a \in \mathcal{A}$, $\delta_{\infty a^\infty}$ is the atomic measure entirely supported on ${}^\infty a^\infty$.

**Bernoulli measures.** Let $v = (v_a)_{a\in\mathcal{A}}$ be a vector of real numbers such that $0 \leq v_a \leq 1$ for all $a \in \mathcal{A}$ and $\sum_{a\in\mathcal{A}} v_a = 1$. Let $\beta_v$ be the discrete probability distribution on $\mathcal{A}$ such that $\beta_v(a) = v_a$ for all $a \in \mathcal{A}$ (a generalisation of the standard Bernoulli law with $n$ outcomes).

The associated *Bernoulli measure* $\mathcal{B}er_v$ on $\mathcal{A}^{\mathbb{Z}}$ is the product measure $\prod_{i\in\mathbb{Z}} \beta_v$, that is,

$$\mathcal{B}er_v([u_0,\ldots,u_n]) = v_{u_0}\cdots v_{u_n} \qquad \text{for all } u_0,\ldots,u_n \in \mathcal{A}^*.$$

In other words, each cell is drawn in an i.i.d. manner according to $\beta_v$. We denote $\mathcal{B}er(\mathcal{A}^{\mathbb{Z}})$ the set of Bernoulli measures on $\mathcal{A}^{\mathbb{Z}}$ *with nonzero* parameters $(v_a)_{a\in\mathcal{A}}$.

**Uniform measure.** In particular, if we take $v_a = \frac{1}{|\mathcal{A}|}$ for all $a \in \mathcal{A}$ in the previous definition, we obtain the *uniform (Bernoulli) measure* $\lambda$.

The *image measure* of $\mu \in \mathcal{M}_\sigma(\mathcal{A}^{\mathbb{Z}})$ by a cellular automaton $(\mathcal{A}^{\mathbb{Z}}, F)$ is defined as $F\mu(B) = \mu(F^{-1}(B))$ for all $B \in \mathfrak{B}$. This defines an action $F : \mathcal{M}_\sigma(\mathcal{A}^{\mathbb{Z}}) \to \mathcal{M}_\sigma(\mathcal{A}^{\mathbb{Z}})$.

We endow $\mathcal{M}_\sigma(\mathcal{A}^{\mathbb{Z}})$ with the *weak-$*$ topology*: for a sequence $(\mu_n)_{n\in\mathbb{N}} \in \mathcal{M}_\sigma(\mathcal{A}^{\mathbb{Z}})^{\mathbb{N}}$ and a measure $\mu \in \mathcal{M}_\sigma(\mathcal{A}^{\mathbb{Z}})$, we have:

$$\mu_n \underset{n\to\infty}{\longrightarrow} \mu \quad \Leftrightarrow \quad \forall u \in \mathcal{A}^*,\ \mu_n([u]) \underset{n\to\infty}{\longrightarrow} \mu([u]).$$

This topology makes $F : \mathcal{M}_\sigma(\mathcal{A}^{\mathbb{Z}}) \to \mathcal{M}_\sigma(\mathcal{A}^{\mathbb{Z}})$ continuous and $\mathcal{M}_\sigma(\mathcal{A}^{\mathbb{Z}})$ compact.

A measure $\mu \in \mathcal{M}_\sigma(\mathcal{A}^{\mathbb{Z}})$ is *$\sigma$-ergodic* if, for every borelian $B \in \mathfrak{B}$ such that $\sigma(B) = B$ $\mu$-almost everywhere, we have $\mu(B) = 0$ or 1. In particular, all examples given earlier are $\sigma$-ergodic and the image of a $\sigma$-ergodic measure under the action of a cellular automaton is $\sigma$-ergodic.

As an example of a non-$\sigma$-ergodic measure, consider the average of two Dirac measures $\frac{1}{2}(\delta_{^\infty 0^\infty} + \delta_{^\infty 1^\infty})$ (the set $\{^\infty 0^\infty\}$ is $\sigma$-invariant and has measure $\frac{1}{2}$).

We make use of the following corollary to Birkhoff's theorem:

Corollary 2. *Let $\mu$ be a $\sigma$-ergodic measure and $u \in \mathcal{A}^*$. Then:*

$$\forall_\mu x \in \mathcal{A}^{\mathbb{Z}},\ \mathrm{freq}(u,x) = \mu([u]),$$

*where $\forall_\mu x$ means for $\mu$-almost all $x$ (that is, for all $x$ in some set of measure* 1*).*

**2. Known and new results.** The first main result on one-dimensional cyclic cellular automata is the following. It describes the evolution the values of the sequence $(C_n^t(x)_0)_{t\in\mathbb{N}}$ for an arbitrary site (here 0) when iterating $C_n$ on a uniform random configuration.

THEOREM 3 (Fisch [10], Theorem 1). *Draw an initial configuration $x$ according to $\lambda$ the uniform Bernoulli measure on $(\mathbb{Z}/n\mathbb{Z})^{\mathbb{Z}}$, and consider the sequence $(C_n^t(x)_0)_{t\in\mathbb{N}}$. Then, $\lambda$-almost surely:*

- *For $n \leq 4$, $C_n^t(x)_0$ changes infinitely often as $t \to \infty$ ($x_0$* fluctuates*);*
- *For $n \geq 5$, $C_n^t(x)_0$ changes finitely often as $t \to \infty$ ($x_0$* fixates*).*

Since changes of values correspond to times when a particle $p_+$ or $p_-$ crosses the central column, this result can be interpreted in terms of limit measures. For $n \geq 5$, some particles $p_=$ ("walls") survive asymptotically ($C_n^t\lambda([p_=]) \not\to 0$) and delimit walled areas where the remaining moving particles $p_-$ or $p_+$ cannot enter; for $n \leq 4$, $C_n^t\lambda([p_=]) \to 0$ and moving particles cross each column infinitely often. This result can be intuited on Figure 2.

Notice that the previous result only applies when the initial measure is uniform. The following result follows from [18], Corollary 1; it is weaker but applies on the much more general setting of $\sigma$-ergodic measures:

PROPOSITION 4. *Let $\mu$ be any $\sigma$-ergodic measure on $(\mathbb{Z}/n\mathbb{Z})^{\mathbb{Z}}$. Then at least two of the following are true:*

- $C_n^t\mu([p_+]) \to 0$*;*
- $C_n^t\mu([p_-]) \to 0$*;*
- $C_n^t\mu([p_=]) \to 0$.

For Bernoulli measures, the state of the art is summed up in the following proposition:

PROPOSITION 5. *If $\mu$ is a Bernoulli measure, then $C_n^t\mu([p_+]) \to 0$ and $C_n^t\mu([p_-]) \to 0$. In particular, if $n = 3$, any limit point of $(C_n^t\mu)_{t\in\mathbb{N}}$ is a convex combination of the measures $\delta_{\infty a^\infty}$, $a \in \mathbb{Z}/n\mathbb{Z}$.*

*If furthermore $\mu = \lambda$ the Bernoulli uniform measure, the unique limit point of $(C_n^t\mu)_{t\in\mathbb{N}}$ is $\frac{1}{n}\sum_{a\in\mathbb{Z}/n\mathbb{Z}}\delta_{\infty a^\infty}$ for both cases $n \in \{3,4\}$.*

PROOF. In the case where $\mu$ is a Bernoulli measure, or more generally a measure invariant by the mirror involution $\gamma : (x_i)_{i\in\mathbb{Z}} \mapsto (x_{-i})_{i\in\mathbb{Z}}$, the only possible non-zero case is $C_n^t\mu([p_=]) \not\to 0$. Indeed, since $C_n \circ \gamma = \gamma \circ C_n$ and the mirror operation sends $p_+$ to $p_-$ and vice versa, we have $C_n^t\mu([p_+]) = C_n^t\mu([p_-])$.

For $n = 3$, since $p_= = \emptyset$, there is asymptotically no particle at all, so all limit points must be some convex combination of the measures $\delta_{\infty a^\infty}$, $a \in \mathbb{Z}/n\mathbb{Z}$.

If furthermore $\mu = \lambda$ the Bernoulli uniform measure, Theorem 3 gives us $C_n^t\mu([p_=]) \to 0$ in the case $n = 4$ as well. Since this measure is invariant by the state-transposing operation $\kappa : (x_i)_{i\in\mathbb{Z}} \mapsto (x_i+1)_{i\in\mathbb{Z}}$ and $C_n \circ \kappa = \kappa \circ C_n$, the unique limit point is $\frac{1}{n}\sum_{a\in\mathbb{Z}/n\mathbb{Z}} \delta_{\infty a^\infty}$ for both cases $n \in \{3, 4\}$. □

The previous results on $C_3$, in particular the statement concerning fluctuation in Theorem 3, can be interpreted in terms of heteroclinic cycles in the orbit of $\lambda$-almost every configuration $x$.

On the one hand, no state ever dominates the whole space in the sense that for every state $a \in \mathbb{Z}/3\mathbb{Z}$,

$$\mathrm{freq}(a, C_3^t(x)) = C_3^t\lambda([a]) \to \frac{1}{3},$$

where the first equation uses Corollary 2 and the fact that the image under $C_3$ of a $\sigma$-ergodic measure is $\sigma$-ergodic, and the second equation uses the last statement of Proposition 5.

On the other hand, if one considers a fixed window $[-N, N]$, Proposition 5 implies that the frequency of non-monochromatic words (i.e. particles) in $C_3^t(x)$ tends to 0. However, Theorem 3 shows that $C_3^t(x)_{[-N,N]}$ never becomes stationary; instead, particles keep crossing the central column, letting each state dominate the window in turn. Still, $C_3^t(x)_{[-N,N]}$ will be monochromatic "most of the time" (in topological terms, it is close to one of the ${}^\infty a^\infty$, $a \in \mathbb{Z}/3\mathbb{Z}$); that is, the increasing sequence $(t_i)_{i\in\mathbb{N}}$ of times where $C_3^{t_i}(x)_{[-N,N]}$ is not monochromatic satisfies $\liminf i/t_i = 0$.

In this sense, the 3-state cyclic cellular automaton exhibits heteroclinic cycles in local regions.

Our main new result determines the unique limit point for non-uniform Bernoulli measures:

THEOREM 6 (Main result).
*Let $\mu$ be a Bernoulli measure on $(\mathbb{Z}/3\mathbb{Z})^\mathbb{Z}$ with nonzero parameters $(p_0, p_1, p_2)$. Then:*

$$C_3^t\mu \underset{t\to\infty}{\longrightarrow} p_2\delta_{\infty 0^\infty} + p_0\delta_{\infty 1^\infty} + p_1\delta_{\infty 2^\infty}.$$

Theorem 6 can be interpreted as follows. Draw an initial configuration according to a Bernoulli measure with nonzero parameters $(p_0, p_1, p_2)$, and consider a fixed arbitrary window $[-N, N]$. By Proposition 5, the probability that $C_3^t(x)_{[-N,N]}$ contains at least two different states (i.e. a particle) tends to 0. Theorem 6 further shows that the probability that $C_3^t(x)_i = a$ for $a \in \mathbb{Z}/3\mathbb{Z}$ and all $i \in [-N, N]$ tends to $p_{a-1}$ as $t$ tends to infinity.

Remarkably, the parameters of the limit measure are a simple cyclic permutation of the parameters of the initial Bernoulli measure: each state $a$ reaches asymptotically the initial frequency of its "prey" $p_{a-1}$. This is illustrated on Figure 3.

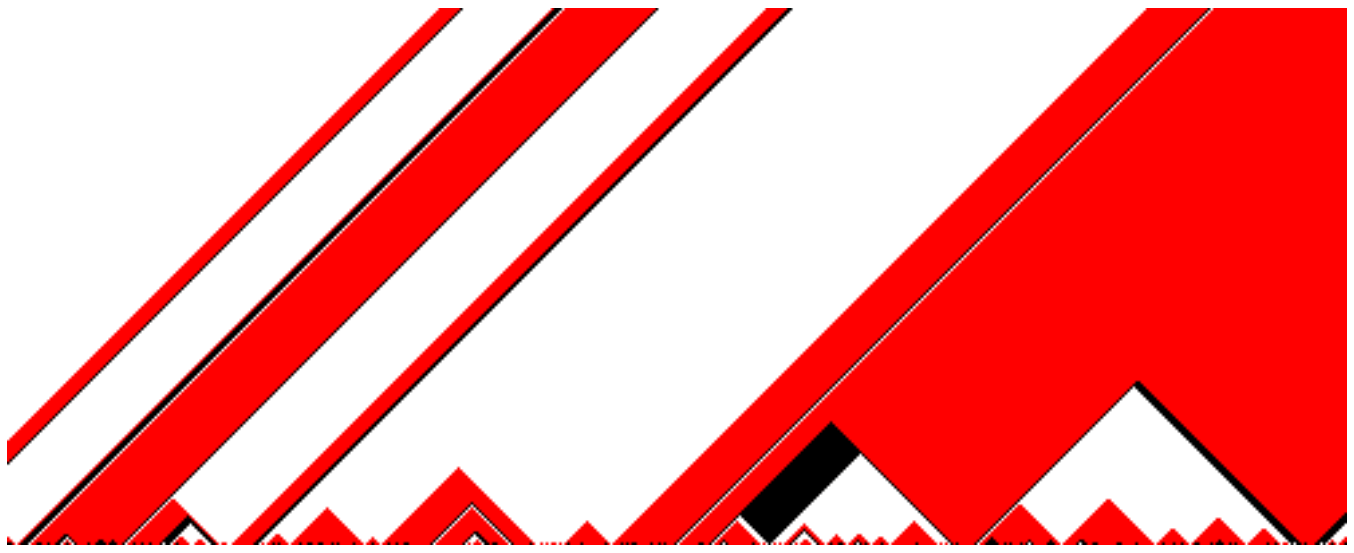

*Fig 3. The 3-state cyclic cellular automaton iterated on an initial configuration drawn according to the Bernoulli measure of parameters $(\frac{1}{10}, \frac{3}{10}, \frac{6}{10})$. State $0 = \square$, initially present with probability $\frac{1}{10}$ at each site, is present in the topmost configuration with a frequency approximately $\frac{6}{10}$.*

## 3. Proof of the main result.

This section is dedicated to the proof of Theorem 6. Since we already know by Proposition 5 that any limit point of $(C_3^t\mu)_{t\in\mathbb{N}}$ is a convex combination of $\delta_{^\infty 0^\infty}$, $\delta_{^\infty 1^\infty}$ and $\delta_{^\infty 2^\infty}$, it remains to show that for each $a$, $\mu(C_3^t(x)_0 = a) \to p_{a-1}$.

In this section, we use the one-sided version of $C_3$ to simplify proofs:

Definition 7 (One-sided cyclic CA). *$(\mathbb{Z}/3\mathbb{Z}, C_{3+})$ is the one-sided 3-state cyclic cellular automaton defined on the neighbourhood $N = \{0, 1\}$ by the local rule:*

$$c_{3+}(u_0, u_1) = \begin{cases} u_0 + 1 & \text{if } u_1 = u_0 + 1 \bmod 3, \\ u_0 & \text{otherwise.} \end{cases}$$

It is easy to check by hand that $c_3(a, b, c) = c_{3+}(c_{3+}(a, b), c_{3+}(b, c))$ for all $a, b, c \in \mathbb{Z}/3\mathbb{Z}$ (assume $a = 0$ by symmetry). Therefore:

$$C_3 = C_{3+}^2 \circ \sigma.$$

Hence proving that $\mu(C_{3+}^t(x)_0 = a) \to p_{a-1}$ implies the same result on $C_3$.

The proof proceeds in 4 steps:

**Section 3.1** where we associate a random walk to each configuration and relate the properties of this random walk to the orbit of the configuration under $C_{3+}$;
**Section 3.2** where we translate Theorem 6 on the random walk and establish the objects that will be relevant to the proof;
**Section 3.3** where we introduce a second random walk "embedded" in the previous one, which is symmetric (hence easier to analyse) and captures its large-scale behaviour;
**Section 3.4** where we bring back the results from the embedded walk to the initial walk and bring all tools together to conclude the proof.

3.1. *Random walk associated with a configuration.* In this section, we introduce tools to turn the study of the orbits of the one-sided 3-state cyclic automaton into the study of some random walk built from the initial configuration $x$.

DEFINITION 8. *To a configuration $x \in \{0,1,2\}^{\mathbb{Z}}$ we associate a random walk $W[x] := (w_i)_{i\in\mathbb{Z}}$ on $\mathbb{Z}$ such that $w_0 \in \{0,1,2\}$ and made up of steps in $\{-1,0,1\}$ as follows:*

- $w_0 = x_0$,
- *for all $i \geq 0$, $w_{i+1}$ is the value in $\{w_i - 1, w_i, w_i + 1\}$ such that $w_{i+1} \equiv x_{i+1} \bmod 3$,*
- *and for $i \leq 0$, $w_{i-1}$ is the value in $\{w_i - 1, w_i, w_i + 1\}$ such that $w_{i-1} \equiv x_{i-1} \bmod 3$.*

*This encoding is an bijection.*

Figure 4 provides an example of this encoding (from the black configuration $x$ to the black walk $W[x]$).

We denote by $W_{[a,b]}[x] := (w_a, w_{a+1}, \ldots, w_b)$ the positions of the walk on $\mathbb{Z}$ from time $a$ to time $b$. Notice that we call time in the context of the random walk what corresponds to space in the configuration $x$, which is different from the time corresponding to the iteration of cellular automaton. Context should make clear which notion of time we refer to.

The main interest of this correspondence is to deduce the state of a cell after $n$ iterations from the maximal height during the first $n$ steps of the walk:

PROPOSITION 9. *For $n \geq 0$, we have*

$$C^n_{3+}(x)_0 = \left(\max W_{[0,n]}[x]\right) \bmod 3.$$

PROOF. We will prove that the following invariant is maintained under the iterations of $C_{3+}$:

$$\left(\max W_{[0,n-t]}[C_{3+}^t(x)]\right) \bmod 3 = \left(\max W_{[0,n-(t+1)]}[C_{3+}^{t+1}(x)]\right) \bmod 3.$$

When this invariant is expressed for $t = 0$ and $t = n$, we deduce the expected identity:

$$\begin{aligned}\left(\max W_{[0,n]}[x]\right) \bmod 3 &= \left(\max W_{[0,0]}[C_{3+}^n(x)]\right) \bmod 3\\ &= \max\{C_{3+}^n(x)_0\} \bmod 3\\ &= C_{3+}^n(x)_0.\end{aligned}$$

We prove this invariant in the case $t = 0$ and any $n \geq 1$. The cases $t > 0$ follow by replacing $x := C_{3+}^t(x)$ and $n := n - t$.

We describe how to obtain $W_{[0,n-1]}[C_{3+}(x)] = (w'_i)_{i=0,\dots,n-1} =: w'$ from $W_{[0,n]}[x] = (w_i)_{i=0,\dots,n} =: w$ by a 3-step transformation

$$w \mapsto w^1 := (w_i^1)_{i=0,\dots,n} \mapsto w^2 := (w_i^2)_{i=0,\dots,n} \mapsto w'.$$

Each of these steps, illustrated in Figure 4, preserves the invariant.

By definition, for $0 \leq i \leq n-1$, $C_{3+}(x)_i = c_{3+}(x_i, x_{i+1})$. We notice that cases where $x_i$ becomes $x_i + 1 \bmod 3$ are exactly the steps $+1$ in the walk $w$ (factors 01, 12 or 20 in $x$).

Step 1 For $0 \leq i < n$, define $w_i^1 := w_i + 1$ if $w_i < w_{i+1}$ and $w_i^1 := w_i$ otherwise. In addition $w_n^1 := w_n$. Notice that $w^1$ is also a walk on $\mathbb{Z}$ made up of steps $\{-1, 0, +1\}$. The maximal height is preserved since any visit at maximal height in $w$ cannot be followed by a $+1$ step.

Step 2 The only case where $w_0^1 = 3 \notin \{0, 1, 2\}$ is when $w_0 = 2 \bmod 3$ and $w_1 = 0 \bmod 3$. In this case, for $i = 0, \dots, n$, define $w_i^2 = w_i^1 - 3$ and $w^2 = w^1$ otherwise. The maximal height may be decreased by 3, but it is preserved mod 3.

Step 3 Remove the last vertex $w_n^2$ from $w^2$ to obtain $w' = W_{[0,n-1]}[C_{3+}(x)]$. This preserves the maximal height: if $w_n$ was the first visit to the maximal height, the first step ensures that $w_{n-1}^1 = w_{n-1} + 1 = w_n = w_n^1$. Therefore $w_{n-1}^2 = w_n^2$, so $w_n^2$ cannot be the first occurrence of the maximal height and can be safely removed.

□

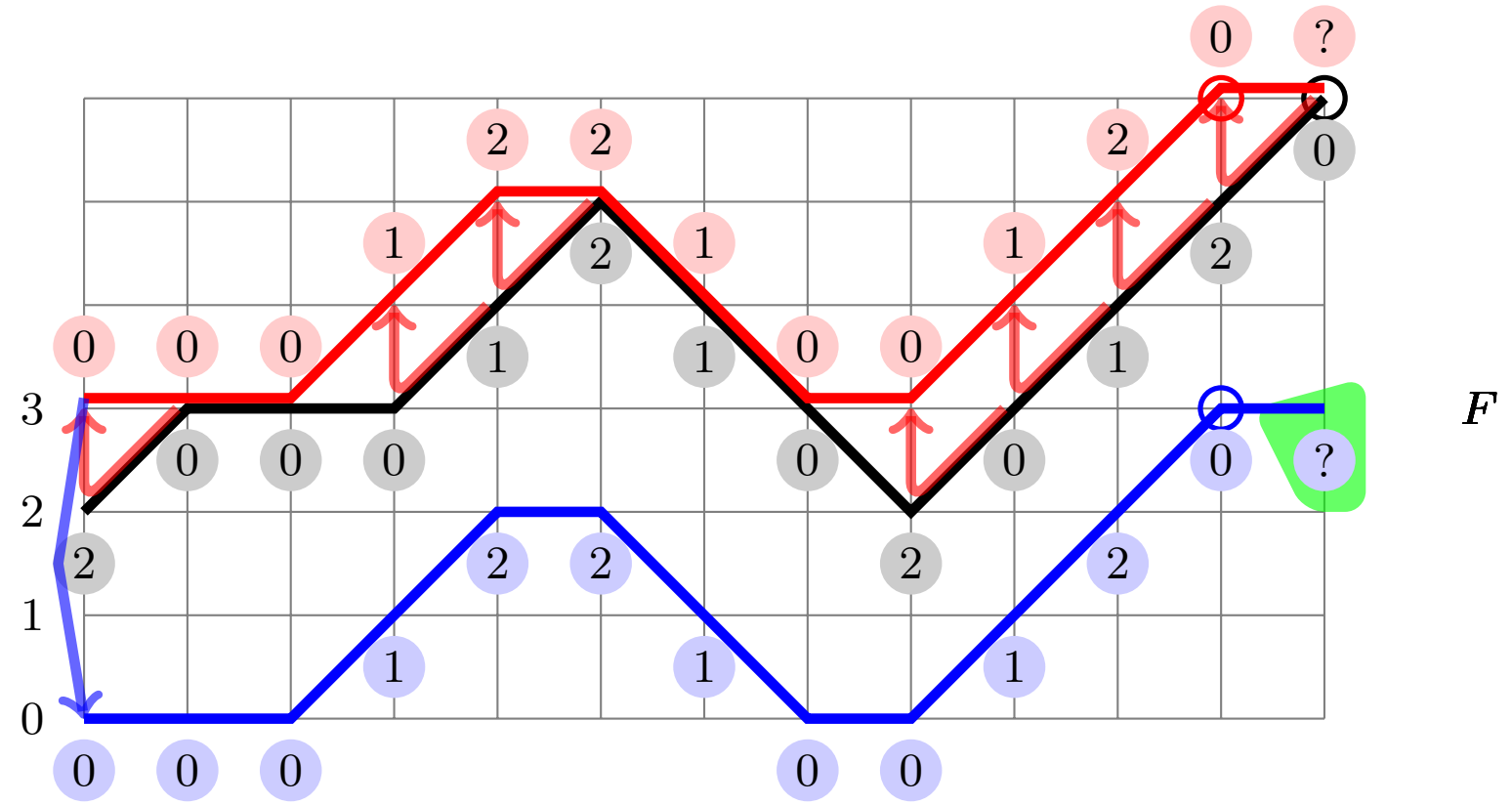

The preserved invariant is the height (modulo 3) of the first maximum of paths represented by ○,○,○, respectively.
The initial configuration given by ($x_i$)$_{i=0,\ldots,12}$ is represented by the black path $(w_i)_{i=0,\ldots,12}$.
Step 1: Its image ($y_i$)$_{i=0,\ldots,12}$ by $C_{3+}$ is represented by the red path $(w_i^1)_{i=0,\ldots,12}$, obtained from the black path by rising in parallel the lower vertex of each rise ( ).
Step 2: Since this red path starts at $3 \notin \{0,1,2\}$, it is shifted downwards by 3, leading to the blue path $(w_i^2)_{i=0,\ldots,12}$ starting at $0 \in \{0,1,2\}$ that is associated with ($y_i$)$_{i=0,\ldots,12}$.
Step 3: The green area indicates that the last vertex of the blue path is deleted, yielding $(w_i')_{i=0,\ldots,11}$.

*FIG 4. The 3-step transformation preserving the invariant.*

3.2. *Analysing the random walk* $\mathbb{P}(Z_t = j) = p_j$. Recall that the measure on the initial configuration is the Bernoulli measure $\mu$ of parameters $(p_0, p_1, p_2)$. From this and the bijection between configuration and walks on $\mathbb{Z}$ we forget its relationship with $x$ to study it for itself as a random variable, directly sampling $W[x]$ as follows (each choice being independent). For $j \in \{0, 1, 2\}$,

- $w_0 = j$ with probability $p_j$;
- for all $i \geq 0$, with probability $p_j$, $w_{i+1}$ is the value in $\{w_i - 1, w_i, w_i + 1\}$ such that $w_{i+1} \equiv j \bmod 3$;
- for $i \leq 0$, with probability $p_j$, $w_{i-1}$ is the value in $\{w_i - 1, w_i, w_i + 1\}$ such that $w_{i-1} \equiv j \bmod 3$.

Similarly, we can sample the factor $W_{[0,n]}[x] = (w_i)_{i=0,\ldots,n}$ by assuming by convention that $w_{-1} = 1$ to ensure that $w_0 \in \{0, 1, 2\}$. Then the only rule is $w_{i+1} \in \{w_i - 1, w_i, w_i + 1\}$ with probability $p_{w_{i+1} \bmod 3}$, independently from other choices.

In the proofs, we need walks that start from an arbitrary $k \in \mathbb{Z}$. Formally, define $\mathcal{W}_{k,n}$ a random walk on $\mathbb{Z}$ of length $n \in \mathbb{N}$ and starting from $k \in \mathbb{Z}$ as $\mathcal{W}_{k,n} := (W_t)_{t=0,\dots,n}$, where:

$$\begin{aligned} &W_0 = k \\ &W_t = W_{t-1} - 1 + ((Z_t - W_{t-1} + 1) \bmod 3) \quad \text{for } t = 1, \dots, n, \end{aligned}$$

where $(Z_t)_{t=1,\dots,n}$ are i.i.d. random variables in $\mathbb{Z}/3\mathbb{Z} := \{0, 1, 2\}$ for all $t$, and $\mathbb{P}(Z_t = j) = p_j$ for all $j \in \mathbb{Z}/3\mathbb{Z}$.

Theorem 10 (Main result of this section). *For any $a \in \mathbb{Z}/3\mathbb{Z}$ and any $k \in \mathbb{Z}$,*

$$\lim_{n\to+\infty} \mathbb{P}(\max(\mathcal{W}_{k,n}) \bmod 3 = a) = p_{(a-1) \bmod 3}$$

*where* $\max(\mathcal{W}_{k,n}) := \max_{t=0,\dots,n} W_t$.

We first consider the case $a = 0$ (and $k = 0$), *i.e.* $\lim_{n\to+\infty} \mathbb{P}(\max(\mathcal{W}_{0,n}) \bmod 3 = 0) = p_2$; the other cases will follow. Our proof proceeds by conditioning this event to the length of the 3-tail (defined below), and describing the probability in terms of the value of other probabilities $(P^{<H}_{k,m})_{k,m,H}$ (also defined below).

Definition 11 (Record, tail). *For a sequence $\mathcal{W} = (W_t)_{0\le t\le n} \in \mathbb{Z}^*$, we say that a* record *occurs at time $t'$ if $W_{t'} = \max_{t=0,\dots,t'} W_t$; notice that a sequence can have multiple records $t_i$ sharing the same height $W_{t_i}$.*

*The $h$-tail of $\mathcal{W}$ is the suffix $W_{[t',n]}$, where $t'$ is the last occurrence of a record whose height $W_{t'}$ is* divisible by *$h$; the $h$-tail for $h > 1$ may not exist. The length of the $h$-tail is denoted by* $\mathsf{tail}_h(\mathcal{W}) := n - t'$.

In this paper we make use of the 3-tail and the 1-tail. $\mathsf{tail}_3(\mathcal{W}_{k,n})$ is usually denoted by $m$.

*Notations.*

- $\mathbb{W}_{k,n}$ is the set of walks on $n \in \mathbb{N}$ steps which start from $k \in \mathbb{Z}$.
- $\mathbb{W}^{<H}_{k,n}$ is the set of walks on $n \in \mathbb{N}$ steps which start from $k \in \mathbb{Z}$ and remain on values strictly lower than $H \in \mathbb{Z}$.
- $\mathbb{W}^{+H}_{k,n}$ is the set of walks in $\mathbb{W}^{<H+1}_{k,n}$ that end on a record of $H$ (not necessarily the first visit).
- $P^{<H}_{k,n}$ and $P^{+H}_{k,n}$ are the probabilities that a random walk $\mathcal{W}_{k,n} \in \mathbb{W}_{k,n}$ belongs to $\mathbb{W}^{<H}_{k,n}$ or $\mathbb{W}^{+H}_{k,n}$, respectively.

Proposition 12 (Description conditioned by 3-tail). *For any $n \in \mathbb{N}$ and any possible 3-tail length $m \in \mathbb{N}$ with $1 \le m \le n$, we have:*

$$\mathbb{P}(\max(\mathcal{W}_{k,n}) \bmod 3 = 0 \mid \mathsf{tail}_3(\mathcal{W}_{k,n}) = m) = p_2 K_m,$$

*where* $K_m := \dfrac{P^{<0}_{-1,m-1}}{P^{<0}_{-1,m}}$.

Proof. By the definition of conditional probability:

$$\mathbb{P}(\max(\mathcal{W}_{k,n}) \bmod 3 = 0 \mid \mathsf{tail}_3(\mathcal{W}_{k,n}) = m) = \frac{\mathbb{P}(\mathsf{tail}_3(\mathcal{W}_{k,n}) = m \text{ and } \max(\mathcal{W}_{k,n}) \bmod 3 = 0)}{\mathbb{P}(\mathsf{tail}_3(\mathcal{W}_{k,n}) = m)}.$$

We now evaluate the denominator and then the numerator of the right-hand side.

Evaluation of $\mathbb{P}(\mathsf{tail}_3(\mathcal{W}_{k,n}) = m)$:
By definition $\mathsf{tail}_3(\mathcal{W}_{k,n}) = m \ge 1$ implies that $n-m$ is the last record whose height $W_{n-m} = 3H$ is divisible by 3 in $\mathcal{W}_{k,n}$, and that $W_{n-m+1}$ exists and is in the 3-tail. By the law of total probability:

$$\mathbb{P}(\mathsf{tail}_3(\mathcal{W}_{k,n}) = m) = \sum_{H \in \mathbb{Z}} P^{+3H}_{k,n-m} \cdot \mathbb{P}(\mathsf{tail}_3(\mathcal{W}_{k,n}) = m \mid \mathcal{W}_{k,n-m} \in \mathbb{W}^{+3H}_{k,n-m}).$$

Assume therefore that $\mathcal{W}_{k,n-m} \in \mathbb{W}^{+3H}_{k,n-m}$ for some $H$, that is, $n-m$ is a record and $W_{n-m} = 3H$. Since $m \ge 1$, we may discuss the possible values of $W_{n-m+1} \in \{3H-1, 3H, 3H+1\}$ for any walk of $\mathcal{W}_{k,n}$. We identify below which of these values are allowed.

- If $W_{n-m+1} = 3H$ (with probability $p_0$), then $n-m$ is not the last record whose value is divisible by 3, a contradiction.
- If $W_{n-m+1} = 3H+1$ (with probability $p_1$), future visits at height $3H$ will not be a new record, so $W_{n-m}$ is the last record divisible by 3 if and only if the walk never reaches $3H+3$. This corresponds to $W_{[n-m+1,n]} \in \mathbb{W}^{<3H+3}_{3H+1,m-1}$ happening with probability $P^{<3H+3}_{3H+1,m-1}$.
- If $W_{n-m+1} = 3H-1$ (with probability $p_2$), the next visit at height $3H$ would be a new occurrence of a record whose value is divisible by 3, so $n-m$ is the last record whose value is divisible by 3 if and only if the walk never visits $3H$ again. This corresponds to $W_{[n-m+1,n]} \in \mathbb{W}^{<3H}_{3H-1,m-1}$ happening with probability $P^{<3H}_{3H-1,m-1}$.

Thus we get:

$$\mathbb{P}(\mathsf{tail}_3(\mathcal{W}_{k,n}) = m) = \sum_{H\in\mathbb{Z}} P^{+3H}_{k,n-m} \cdot (p_1 P^{<3H+3}_{3H+1,m-1} + p_2 P^{<3H}_{3H-1,m-1}).$$

In the definition of $\mathcal{W}_{k,n}$ it appears that any realisation $(W_t)_t \in \mathcal{W}_{k,n}$ can be translated into $(W_t+3T)_t$ for any $T \in \mathbb{Z}$ without changing the probability of steps. This implies that for any $3T \in 3\mathbb{Z}$ we have:

$$\forall (H,k) \in \mathbb{Z}^2, \forall n \in \mathbb{N}, P^{<H}_{k,n} = P^{<H+3T}_{k+3T,n}.$$

Therefore the probabilities in the previous discussion do not depend on $3H$.

$$\mathbb{P}(\mathsf{tail}_3(\mathcal{W}_{k,n}) = m) = \left(\sum_{H\in\mathbb{Z}} P^{+3H}_{k,n-m}\right) \cdot (p_1 P^{<0}_{-2,m-1} + p_2 P^{<0}_{-1,m-1})$$

The first step of a walk in $\mathbb{W}^{<0}_{-1,m}$ goes from $-1$ to either $-1$ or $-2$, so we get the following partition:

$$\mathbb{W}^{<0}_{-1,m} = \{-1\} \times \mathbb{W}^{<0}_{-1,m-1} \cup \{-1\} \times \mathbb{W}^{<0}_{-2,m-1}$$

In terms of probabilities this identity turns into:

$$P^{<0}_{-1,m} = p_2 P^{<0}_{-1,m-1} + p_1 P^{<0}_{-2,m-1}.$$

Hence $\mathbb{P}(\mathsf{tail}_3(\mathcal{W}_{k,n}) = m) = \left(\sum_{H\in\mathbb{Z}} P^{+3H}_{k,n-m}\right) P^{<0}_{-1,m}$.

Evaluation of $\mathbb{P}(\mathsf{tail}_3(\mathcal{W}_{k,n}) = m \text{ and } \max(\mathcal{W}_{k,n}) \bmod 3 = 0)$: We adapt the previous discussion on $W_{n-m+1}$ under the additional condition $\max(\mathcal{W}_{k,n}) \bmod 3 = 0$. Again assume that $\mathcal{W}_{k,n-m} \in \mathbb{W}^{+3H}_{k,n-m}$ for some $H$.

- $W_{m-n+1} = 3H$ is still impossible by definition of the 3-tail.
- $W_{m-n+1} = 3H+1$ is impossible since it would imply the maximum is strictly greater that $3H$. However, there can be no new record whose value is divisible by 3 by definition of 3-tail.

Hence the only possible choice for $W_{n-m+1}$ is $3H-1$, and the walk must avoid $3H$ from time $n-m+1$ onwards: this happens with probability $p_2 P^{<3H}_{3H-1,m-1} = p_2 P^{<0}_{-1,m-1}$. Therefore the numerator of $K_m$ is:

$$\mathbb{P}(\mathsf{tail}_3(\mathcal{W}_{k,n}) = m \text{ and } \max(\mathcal{W}_{k,n}) \bmod 3 = 0) = \left(\sum_{H\in\mathbb{Z}} P^{+3H}_{k,n-m}\right) p_2 P^{<0}_{-1,m-1}$$

It is clear that $\left(\sum_{H\in\mathbb{Z}} P^{+3H}_{k,n-m}\right) \neq 0$, so they cancel out in the expression for $K_m$ to give the desired result. □

3.3. *The embedded walk.* In the remainder of this section, we prove that when $n$ tends to infinity, the length $m$ of the 3-tail also tends to infinity with high probability (Lemma 17) and $K_m$ tends to 1 (Lemma 18). Proposition 12 then leads to Theorem 10 for $i = 0$.

We use a factorisation of the walk into a symmetric $\{+3, -3\}$ random walk on $3\mathbb{Z}$, that can be scaled to be the usual symmetric $\{+1, -1\}$ random walk on $\mathbb{Z}$.

DEFINITION 13 (Embedded walk). *Define inductively a sequence of times $(t_j)_{j=0,\dots,\ell-1}$ as follows:*

- *$W_{t_0}$ is the first occurrence of a value divisible by 3 in $\mathcal{W}_{k,n}$, if any.*
- *Given $(t_j)_{j=0,\dots,i-1}$, $t_i > t_{i-1}$ is the next time when $W_{t_i}$ is divisible by 3 and distinct from $W_{t_{i-1}}$, if any such $t_i \le n$ exists.*

*From this sequence we define the embedded walk map $\mathsf{emb} : \mathbb{W}_{k,n} \to (3\mathbb{Z})^*, (W_t)_t \mapsto (W_{t_j})_{j=0,\dots J-1}$ and the embedded random walk $\mathcal{W}^s_{k,n} = \mathsf{emb}(\mathcal{W}_{k,n})$ whose length $|\mathcal{W}^s_{k,n}|$, also denoted $\ell$, is random.*

An example of embedded walk is given in Figure 5.

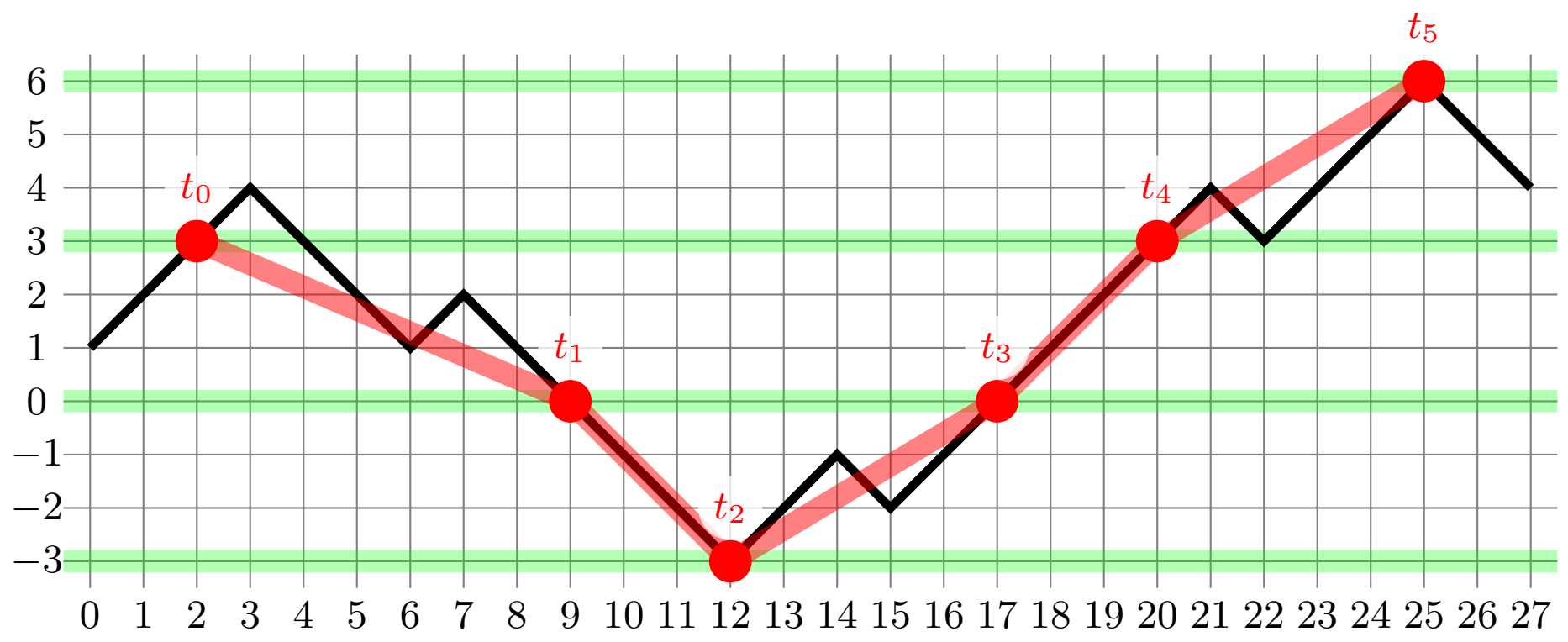

FIG 5. *A (black) walk $\mathcal{W}_{1,27}$ of length 27 and its embedded (red) walk $\mathcal{W}^s_{1,27}$ of length $\ell = 6$.*

Any $\{+3, -3\}$-embedded walk on $3\mathbb{Z}$ corresponds by a scaling of $1/3$ to a $\{+1, -1\}$-walk on $\mathbb{Z}$. In combinatorics such walks are described by words on the alphabet $\{a, b\}$ where the letter $a$ denotes a $+1$ step and the letter $b$ a $-1$ step. Denote $\mathsf{Bin}_\ell := \{a, b\}^\ell$ the set of binary words of length $\ell$ (unconstrained walks). For example, the embedded (red) walk in Figure 5 corresponds to the word *bbaaa*.

First we show that the embedded walk on $3\mathbb{Z}$ is symmetric with independent steps (Lemma 14), then that its length is at least linear in the size of the original walk with high probability (Lemma 15) and finally that the probability of a small 1-tail in a symmetric $\{+1,-1\}$ walk on $\mathbb{Z}$ is asymptotically negligible (Lemma 16).

LEMMA 14 (The embedded walk on $3\mathbb{Z}$ is symmetric). *For any integers $n \geq 0$ and $0 < \ell \leq \lfloor n/3 \rfloor$, the walk $\mathcal{W}^s_{k,n}$ conditioned by $|\mathcal{W}^s_{k,n}| = \ell$ is a symmetric random walk of $\ell - 1$ steps $\{-3,+3\}$ on $3\mathbb{Z}$ with independent steps.*

Notice that the following proof contains the definition of $\tau$-equivalence that will be later used in the proof of Lemma 18.

PROOF. The vertex $W_{t_0}$ exists since the length of the embedded walk is conditioned to be $\ell > 0$. We prove that $\mathcal{W}^s_{k,n}$ is symmetric, that is, for any $i$ such that $1 \leq i \leq \ell - 1$:

$$\mathbb{P}(W_{t_i} = W_{t_{i-1}} + 3) = \mathbb{P}(W_{t_i} = W_{t_{i-1}} - 3)$$

The proof of this equality relies on a set of involutions $\tau_i : (\mathbb{Z}/3\mathbb{Z})^* \to (\mathbb{Z}/3\mathbb{Z})^*$ acting on the variables $(Z_k)_{k=t_{i-1}+1,\ldots,t_i} =: Z_{]t_{i-1},t_i]}$. We illustrate this involution on Figure 6.

First, for any sequence $z = (Z_i)_{0\leq i\leq \ell-1} \in \mathbb{Z}/2\mathbb{Z}^*$, define the mirror image of $z$ as $\overline{z} := (Z_{\ell-i-1})_{0\leq i\leq \ell-1}$.

Second, let $(W_i)_{0\leq i\leq \ell-1}$ be the walk associated to $z$ and define the sequence of times $(t_j)$ as in Definition 13. If $t_j$ is defined, let $s_j$ be the index of the last occurrence of the event $Z_t = 0$ (that is, $W_t = 0 \bmod 3$) before $t_j$; notice $s_j \geq t_{j-1}$ since $Z_{t_{j-1}} = 0$. In that case, put:

$$\tau_j : Z_{[0,\ell]} \mapsto Z_{[0,s_j]}\overline{Z_{]s_j,t_j[}}Z_{[t_j,\ell]},$$

and leave $\tau_j(z)$ undefined if $t_j$ is undefined.

$\tau_j$ is a fixpoint-free involution from the events where $W_{t_j} = W_{t_{j-1}} + 3$ to the events where $W_{t_j} = W_{t_{j-1}} - 3$ (as illustrated in Figure 6). Notice that all maps $\tau_j$ commute, when they are defined, since they depend on and modify different intervals, and do not impact the values of $(t_j)$. We say that two walks $w$ and $w'$ are $\tau$-equivalent if and only if there exists a subset $T$ of the times $(t_j)$ such that $w = \left(\prod_{i\in T}\tau_i\right)(w')$ where the product denotes composition. This is an equivalence relation since the $\tau_j$ are commuting involutions.

Denote $C_\tau(w)$ the $\tau$-equivalence class of a word $w$. Notice that:

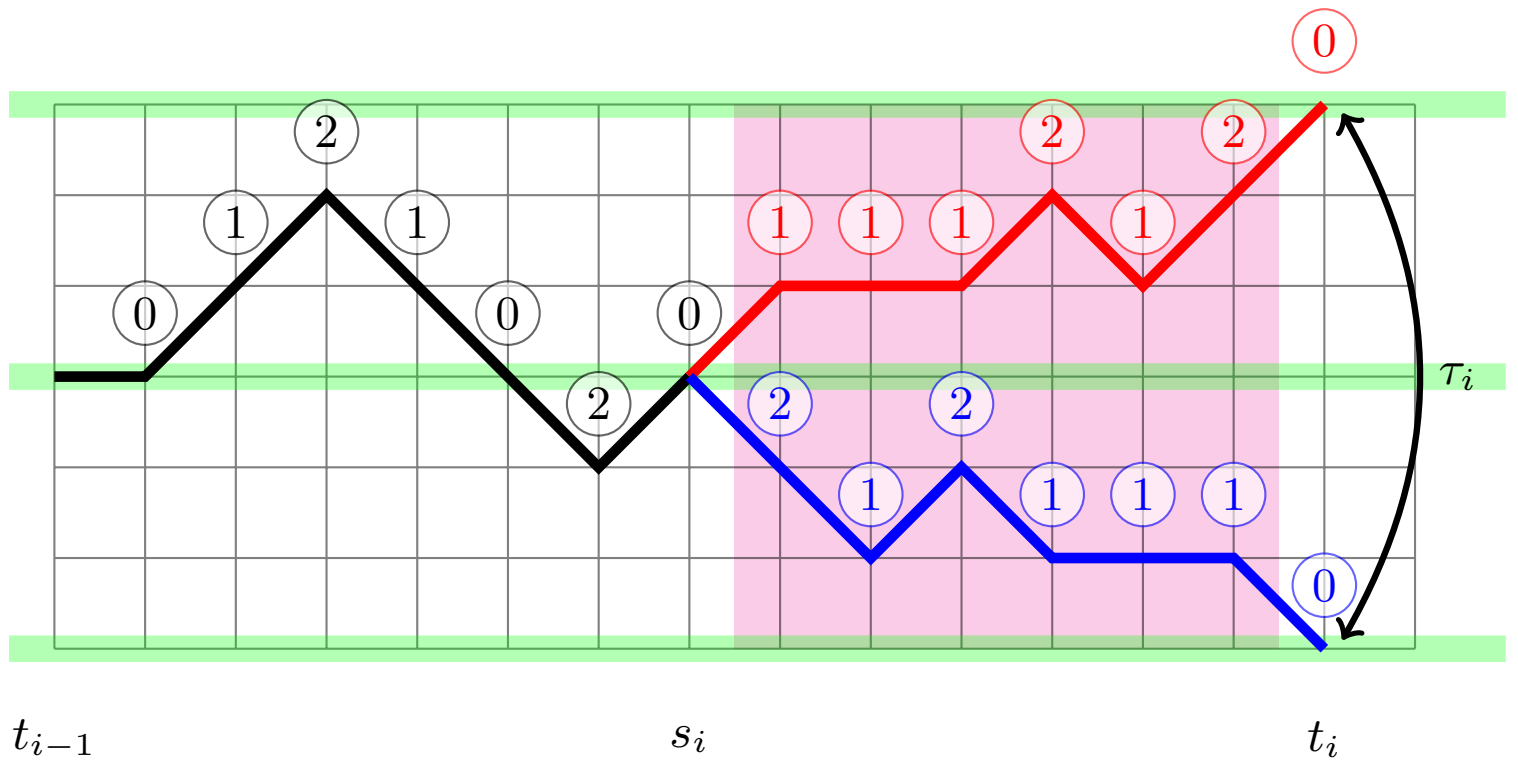

Consider an initial walk (on top) defined by $(Z_t)_{t=t_{i-1}+1,\ldots,t_i} = (\textcircled{x})_t \cup (\textcircled{x})_t$ from $t_{i-1}$ to $t_i$. It defines a $+3$ step in the embedded walk. The magenta interval starts at time $s_i$, the last visit at height 0, and ends at time $t_i$. In this interval, $\tau_i$ transforms the (top) $(\textcircled{x})_t$ into the (bottom) $(\textcircled{x})_t$.
Oberve that the blue and red parts are, up to translation, symmetric with respect to a vertical axis and now defines a step from 0 to $-3$ in the embedded walk.

*Fig 6. Involution used in the proof of symmetry of the embedded walk*

1. Since $Z_{]a,b]}$ and $\overline{Z}_{]a,b]}$ have the same probability, $\tau_j$ preserves probabilities. Therefore, each walk in $C_\tau(w)$ has the same probability.
2. Since $\tau_j$ are fixpoint-free commuting involutions, every such class contains $2^{\ell-1}$ walks where $\ell$ is the length of the embedded walk for all walks.
3. Each word in $\mathsf{Bin}_{\ell-1}$ is represented by exactly one walk in $C_\tau(w)$.

Hence, conditioned to belong to the equivalent class $C_\tau(w)$, the embedded walks are all equiprobable. The sum over all $\tau$-equivalent classes preserves this equiprobability of embedded walks. Thus we proved that the embedded walk has independent steps. □

The following lemma ensures that the length of the embedded walk grows almost surely at least linearly in $n$. This helps us later to convert bounds on 1-tail length in $\frac{1}{3}\mathcal{W}^s_{k,n}$ into bounds on 3-tail length in $\mathcal{W}_{k,n}$.

Lemma 15 (The embedded walk's length is almost always linear). *For any $\beta$ such that $0 < \beta < p_0^2 p_1 p_2/2$, there exists $N \in \mathbb{N}$ such that for any $n \geq \mathbb{N}$ we have:*

$$\mathbb{P}(|\mathcal{W}^s_{k,n}| < \beta n) < \exp(-C(\beta)n),$$

*where $C(\beta) > 0$ is a non-negative function of $\beta$ made explicit in the proof.*

Proof. Split the variables $(Z_i)_{i=1,\dots,n}$ used in the definition of the walk into independent factors of 4 variables $(Z_{]4a,4(a+1)]})_{a=0,\dots,\lfloor n/4\rfloor}$. Now define variables $(X_a)_{a=0,\dots,\lfloor n/4\rfloor}$ as follows:

$$X_a = \begin{cases} 1 & \text{if } Z_{]4a,4(a+1)]} \in \{(0,1,2,0),(0,2,1,0)\} \\ 0 & \text{otherwise} \end{cases} .$$

The motivation for this definition is that, as illustrated in Figure 7, each occurrence of a factor $Z_{]a,a+4]} = (0,1,2,0)$ or $Z_{]a,a+4]} = (0,2,1,0)$ implies the existence of at least one embedded vertex at $t_i = a+4$. This is obviously a rough bound and many embedded vertices are missed.

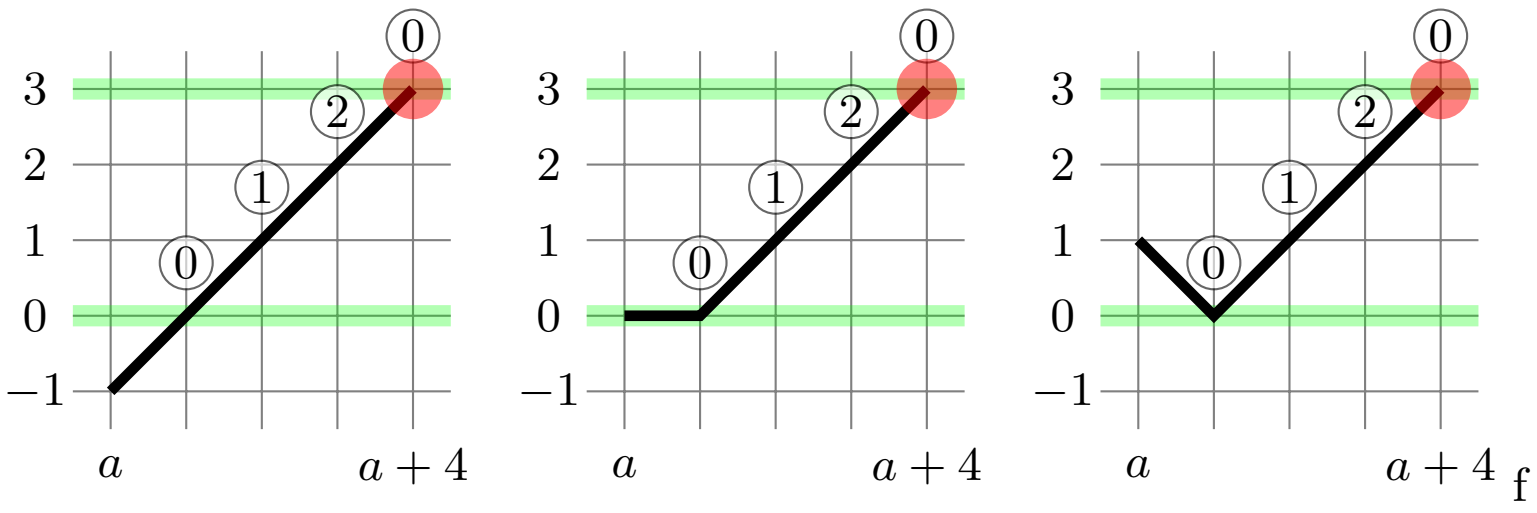

*Fig 7. If $Z_{]a,a+4]} = (0,1,2,0)$ then an embedded step occurs at $t_i = a+4$ (and possibly another at $t_{i-1} = a+1$).*

It follows that $|\mathcal{W}^s_{k,n}| \geq \sum_{a=0}^{\lfloor n/4\rfloor} X_a$. Notice that $(X_a)_{a=0,\dots,\lfloor n/4\rfloor}$ is a family of i.i.d Bernoulli variables of parameter $q := 2p_0^2p_1p_2 < 1$, so $\sum_{a=0}^{\lfloor n/4\rfloor} X_a$ has expectation $\lfloor n/4\rfloor q$. Since $4\beta < q$, Take $N$ large enough that $\frac{4\beta}{1-4/N} \leq \frac{4\beta+q}{2}$. For $n \geq N$, we have:

$$\mathbb{P}\left(\sum_{a=0}^{\lfloor n/4\rfloor} X_a \leq \beta n\right) \leq \mathbb{P}\left(\sum_{a=0}^{\lfloor n/4\rfloor} X_a \leq \frac{4\beta+q}{2}\left\lfloor\frac{n}{4}\right\rfloor\right).$$

By Hoeffding's inequality (see for example [4], Chapter 2.6),

$$\mathbb{P}\left(\sum_{a=0}^{\lfloor n/4\rfloor} X_a \leq \frac{4\beta+q}{2}\left\lfloor\frac{n}{4}\right\rfloor\right) \leq \exp\left(-2\left(\frac{q-4\beta}{2}\right)^2\left\lfloor\frac{n}{4}\right\rfloor\right).$$

Since $|\mathcal{W}^s_{k,n}| \geq \sum_{a=0}^{\lfloor n/4\rfloor} X_a$, this is the desired result. □

To show that the 3-tail of $\mathcal{W}_{k,n}$ is almost always sufficiently large, we consider the 1-tail of the symmetric embedded walk $\mathcal{W}^s_{k,n}$ of length $|\mathcal{W}^s_{k,n}|$. Denote $\mathcal{S}^n = (S^n_t)_{t=0,\dots,n}$ the usual symmetric random walk on $\mathbb{Z}$ made up of $n$ steps $\{-1,+1\}$ (where $S^n_0 = 0$).

LEMMA 16 (Upper bound on the probability of a small 1-tail). *For any* $\alpha \in ]0,1[$,

$$\mathbb{P}(\mathsf{tail}_1(\mathcal{S}^n) \leq n^\alpha) \xrightarrow[n\to\infty]{} 0.$$

PROOF. The limit law of $\mathsf{tail}_1(\mathcal{S}^n)$ as $n \to \infty$ is given by the third arcsine law for the Wiener process; see [8], Chapter III.4 or [12], Section IX.11.1 for more details. For any $0 \leq K \leq 1$,

$$\mathbb{P}(\mathsf{tail}_1(\mathcal{S}^n) \leq K \cdot n) \xrightarrow[n\to\infty]{} \frac{2}{\pi} \arcsin \sqrt{K}.$$

Choose some $\varepsilon > 0$. There exists $\delta > 0$ such that for any $K \leq \delta$, $\frac{2}{\pi} \arcsin \sqrt{K} \leq \frac{\varepsilon}{2}$. Then for $n$ large enough, we have $n^\alpha \leq \delta \cdot n$, so:

$$\begin{aligned}\mathbb{P}(\mathsf{tail}_1(\mathcal{S}^n) \leq n^\alpha) &\leq \mathbb{P}(\mathsf{tail}_1(\mathcal{S}^n) \leq K \cdot n) \\ &\leq 2\pi \arcsin \sqrt{K} + \frac{\varepsilon}{2} \leq \varepsilon.\end{aligned}$$

This is true for any $\varepsilon > 0$, so the result follows. □

3.4. *Back to the main walk, and end of the proof.* We now transfer the bound on the probability of a small 1-tail for the symmetric embedded walk, obtained in Lemma 16, to a similar upper bound for the 3-tail on the initial walk, using the probabilistic lower bound on the length of the embedded walk obtained in Lemma 15.

LEMMA 17 (Upper bound for the probability of a small 3-tail). *For any* $\alpha \in ]0,1[$ *and any* $\beta \in ]0, p_0^2 p_1 p_2/2[$, *we have*

$$\lim_{n\to+\infty} \mathbb{P}(\mathsf{tail}_3(\mathcal{W}_{k,n}) \leq (\beta n)^\alpha) = 0.$$

PROOF. We discuss the probability that $\mathsf{tail}_3(\mathcal{W}_{k,n}) \leq (\beta n)^\alpha$ in $\mathcal{W}_{k,n}$ by conditioning on the length of the embedded walk $\mathcal{W}^s_{k,n}$:

$$\begin{aligned}\mathbb{P}(\mathsf{tail}_3(\mathcal{W}_{k,n}) \leq (\beta n)^\alpha) &= \mathbb{P}(\mathsf{tail}_3(\mathcal{W}_{k,n}) \leq (\beta n)^\alpha \mid |\mathcal{W}^s_{k,n}| < \beta n) \cdot \mathbb{P}(|\mathcal{W}^s_{k,n}| < \beta n) \\ &+ \mathbb{P}(\mathsf{tail}_3(\mathcal{W}_{k,n}) \leq (\beta n)^\alpha \mid |\mathcal{W}^s_{k,n}| \geq \beta n) \cdot \mathbb{P}(|\mathcal{W}^s_{k,n}| \geq \beta n).\end{aligned}$$

For the first term, according to Lemma 15, $\mathbb{P}(|\mathcal{W}^s_{k,n}| < \beta n) \leq \exp(-C(\beta)n)$ which tends to 0 when $n$ tends to infinity.

For the second term, notice that if $\mathsf{tail}_3(\mathcal{W}_{k,n}) \leq (\beta n)^\alpha$ then in particular $\mathsf{tail}_1(\mathcal{W}^s_{k,n}) \leq \frac{1}{3}(\beta n)^\alpha$. It follows that:

$$\begin{aligned}\mathbb{P}(\mathsf{tail}_3(\mathcal{W}_{k,n}) \leq (\beta n)^\alpha \mid |\mathcal{W}^s_{k,n}| \geq \beta n) &\leq \sum_{\ell=\beta n}^{n/3} \mathbb{P}(\mathsf{tail}_1(S^\ell) \leq \frac{1}{3}(\beta n)^\alpha) \cdot \mathbb{P}(|\mathcal{W}^s_{k,n}| = \ell) \\ &\leq \max_{\ell \in [\beta n, n/3]} \mathbb{P}(\mathsf{tail}_1(S^\ell) \leq \ell^\alpha),\end{aligned}$$

where we used the fact that $\mathcal{W}^s_{k,n}$ conditioned on $|\mathcal{W}^s_{k,n}| = \ell$ is a usual symmetric random walk of $\ell$ steps (by Lemma 14) and that $(\beta n)^\alpha \leq \ell^\alpha$.

According to Lemma 16, $\lim_{\ell\to+\infty} \mathbb{P}(\mathsf{tail}_1(S^\ell) \leq \ell^\alpha) = 0$ so

$$\lim_{n\to+\infty} \max_{l\in[\beta n, n/3]} \mathbb{P}(\mathsf{tail}_1(S^\ell) \leq \ell^\alpha) = 0.$$

As both terms tend to 0, we have proved the result. □

Recall that

$$K_m := \frac{P^{<0}_{-1,m-1}}{P^{<0}_{-1,m}},$$

where $P^{<H}_{k,n}$ is the probability that $\mathcal{W}_{k,n}$ remains strictly below $H$. In the following, we establish bounds on $K_m$ by characterising this event in terms of the binary word associated with $\mathsf{emb}(w)$.

We denote by $|w|_x$ the number of occurrences of the letter $x$ in the word $w$ and by $\mathsf{flip}$ the flip map that replaces each occurence of a letter $a$ by a letter $b$ and conversely: e.g. $\mathsf{flip}(abba) = baab$. The proof will involve the following sets of words:

- the set $\mathsf{Dyck}_\ell$ of Dyck words with $\ell$ letters: the word $w$ is a Dyck word if $|w|_a = |w|_b$ and for any prefix $p$, $|p|_a \geq |p|_b$ (excursions).
- the set $\mathsf{Pref}_\ell$ of prefixes of Dyck words with $\ell$ letters (meanders).

Cardinalities of those sets are known: $|\mathsf{Bin}_\ell| = 2^\ell$, $|\mathsf{Dyck}_{2\ell}| = \frac{1}{\ell+1}\binom{2\ell}{\ell}$, $|\mathsf{Dyck}_{2\ell+1}| = 0$, $|\mathsf{Pref}_{2\ell}| = \binom{2\ell}{\ell}$ and $|\mathsf{Pref}_{2\ell+1}| = 2\binom{2\ell}{\ell} - \frac{1}{\ell+1}\binom{2\ell}{\ell}$ (see [12], Section 1.5.3 for references). We will use the upper bound $|\mathsf{Dyck}_\ell| \leq \frac{1}{\ell+1}|\mathsf{Pref}_\ell|$ for any $\ell \geq 0$.

We also use the notation $f(n) = \Omega(g(n))$ to mean that there exists a constant $C > 0$ such that $f(n) \geq Cg(n)$ for all $n \in \mathbb{N}$.

LEMMA 18 (Bounds for $K_m$). *There exist positive constants $\gamma, M$ such that*

$$\forall m \geq M,\ 1 \leq K_m \leq 1 + \frac{\gamma}{m}.$$

PROOF. First, the bound $1 \leq K_m$ is obvious by inclusion of events.

For a walk $\mathcal{W}_{-1,m} = (w_i)_{i=0,\dots m-1}$, denote $\mathcal{W}^-_{-1,m}$ the prefix $w_{[0,m-2]}$ excluding only the last vertex and $\mathcal{W}^{-s}_{-1,m}$ the embedded walk associated with the prefix.

We condition the probabilities $P^{<0}_{-1,m-1}$ and $P^{<0}_{-1,m}$ by the length of their embedded walk; for the latter, we remove the last vertex so that the embedded walks have the same distribution.

$$P^{<0,\ell}_{-1,m-1} := \mathbb{P}(\mathcal{W}_{-1,m-1} \in \mathbb{W}^{<0}_{-1,m-1} \mid |\mathcal{W}^s_{-1,m-1}| = \ell);$$

$$P^{<0,\ell}_{-1,m} := \mathbb{P}(\mathcal{W}_{-1,m} \in \mathbb{W}^{<0}_{-1,m} \mid |\mathcal{W}^{-s}_{-1,m}| = \ell).$$

Since $\mathcal{W}^{-}_{-1,m}$ is distributed as $\mathcal{W}_{-1,m-1}$ by definition, we have:

$$P^{<0}_{-1,m-1} - P^{<0}_{-1,m} = \sum_{\ell} \mathbb{P}(|\mathcal{W}^{s}_{-1,m-1}| = \ell) \left(P^{<0,\ell}_{-1,m-1} - P^{<0,\ell}_{-1,m}\right).$$

Let $\beta \in ]0, p_0^2 p_1 p_2/2[$ and $\ell_{cut} := \lceil \beta(m-1) \rceil$. We split the previous sum at $\ell_{cut}$. For indices $\ell < \ell_{cut}$, the sum is upper bounded by $\mathbb{P}(\mathcal{W}^{s}_{-1,m-1} < \ell_{cut}) \leq \exp(-C(\beta)m)$ according to Lemma 15.

For indices $\ell \geq \ell_{cut}$, we consider each term $P^{<0,\ell}_{-1,m-1} - P^{<0,\ell}_{-1,m}$ individually. This term is the probability that $\mathcal{W}_{-1,m}$ satisfies $W_{m-1} = 0$ assuming $\mathcal{W}^{-}_{-1,m} \in \mathbb{W}^{<0}_{-1,m-1}$ and $|\mathcal{W}^{-s}_{-1,m}| = \ell$. This implies that $W_{t_0} = -3 = W_{t_{\ell-1}}$ and $W_{t_i} \leq -3$ for $i = 1, \ldots \ell - 2$.

In terms of the word describing the embedded walk, this corresponds to $\mathcal{W}^{-s}_{-1,m} \in \mathsf{flip}(\mathsf{Dyck}_{\ell-1})$, so:

$$P^{<0,\ell}_{-1,m-1} - P^{<0,\ell}_{-1,m} \leq \mathbb{P}(W_{t_0} = -3 \text{ and } \mathcal{W}^{-s}_{-1,m} \in \mathsf{flip}(\mathsf{Dyck}_{\ell-1}) \mid |\mathcal{W}^{-s}_{-1,m}| = \ell).$$

For similar reasons, $\mathcal{W}^{-}_{-1,m} \in \mathbb{W}^{<0}_{-1,m-1}$ if and only if $W_{t_0} = -3$ and $\mathcal{W}^{-s}_{-1,m} \in \mathsf{flip}(\mathsf{Pref}_{\ell-1})$. Hence,

$$P^{<0,\ell}_{-1,m-1} = \mathbb{P}(W_{t_0} = -3 \text{ and } \mathcal{W}^{-s}_{-1,m} \in \mathsf{flip}(\mathsf{Pref}_{\ell-1}) \mid |\mathcal{W}^{-s}_{-1,m}| = \ell).$$

Now, since all walks in a $\tau$-equivalence class have the same value for $W_{t_0}$, the event $W_{t_0} = -3$ is independent from any event related to $\mathcal{W}^{-s}_{-1,m}$ under the conditioning $|\mathcal{W}^{-s}_{-1,m}| = \ell$. By Lemma 14, the two previous equations, and the upper bound on the ratio between the number of Dyck words and the number of Dyck prefixes,

$$\frac{P^{<0,\ell}_{-1,m-1} - P^{<0,\ell}_{-1,m}}{P^{<0,\ell}_{-1,m-1}} \leq \frac{|\mathsf{Dyck}_{\ell-1}|}{|\mathsf{Pref}_{\ell-1}|} \leq \frac{1}{\ell}.$$

Then in the sum of terms for $\ell \geq \ell_{cut}$,

$$\begin{aligned}\sum_{\ell \geq \ell_{cut}} \mathbb{P}(\ell(\mathcal{W}_{-1,m-1}) = \ell) \left(P^{<0,\ell}_{-1,m-1} - P^{<0,\ell}_{-1,m}\right) \\ \leq \sum_{\ell \geq 0} \mathbb{P}(\ell(\mathcal{W}_{-1,m-1}) = \ell) \left(\frac{P^{<0,\ell}_{-1,m-1}}{\ell_{cut}}\right) \\ \leq \frac{1}{\ell_{cut}} P^{<0}_{-1,m-1}.\end{aligned}$$

Gathering the bounds on both parts of the sum, we have:

$$P^{<0}_{-1,m-1} - P^{<0}_{-1,m} \leq \exp(-C(\beta)m) + \frac{1}{\beta m} P^{<0}_{-1,m-1}.$$

We now show that the term $\exp(-C(\beta)m)$ is asymptotically negligible compared to $\frac{1}{\beta m} P^{<0}_{-1,m-1}$.

Take the set of walks $w$ such that $w_2 = -3$ (this happens with probability $p_2 p_0$) and whose embedded walk $\mathsf{emb}(w)$ belongs to $\mathsf{Pref}_{\ell-1}$; these walks belong to $\mathbb{W}^{<0}_{-1,m-1}$. By equiprobability in the $\tau$-class of $w$, $\mathsf{emb}(w)$ belongs in $\mathsf{Pref}_{\ell-1}$ with probability $|\mathsf{Pref}_{\ell-1}|/|\mathsf{Bin}_{\ell-1}| = \binom{2(\ell-1)}{\ell-1}/2^{\ell-1}$. The minimum value for this probability is when $\ell$ takes its maximum possible value $\ell = m/3$, so its asymptotic behaviour is of order $\sqrt{3/m}$ at least. Therefore $\frac{1}{\beta m} P^{<0}_{-1,m-1}$ is at least of order $\Omega(m^{-3/2})$.

Hence asymptotically, there exists $M$ and $\alpha > 1$ such that for any $m \geq M$,

$$P^{<0}_{-1,m-1} - P^{<0}_{-1,m} \leq \frac{\alpha}{\beta m} P^{<0}_{-1,m-1}.$$

Dividing this equation by $P^{<0}_{-1,m}$, we obtained the desired upper bound on $K_m$:

$$K_m \leq \frac{1}{1 - \frac{\alpha}{\beta m}} \leq 1 + \frac{\gamma}{m}$$

for any $\gamma > \alpha/\beta$ and $m$ large enough. □

Proof. (Theorem 6 for $a = 0$ and $k = 0$). Proposition 12 describes the probability conditioned by the length of the 3-tail:

$$\mathbb{P}(\max \mathcal{W}_{k,n} \bmod 3 = 0) = \sum_{m \geq 0} p_2 K_m \mathbb{P}(\mathsf{tail}_3(\mathcal{W}_{k,n}) = m).$$

By Lemma 18, $K_m \geq 1$ so $\mathbb{P}(\max \mathcal{W}_{k,n} \bmod 3 = 0) - p_2 \geq 0$.

Let $\alpha \in ]0,1[$ and $\beta \in ]0, p_0^2 p_1 p_2/2[$. We split the previous sum around the threshold $m_{cut} := (\beta n)^\alpha$.

- For $m \leq m_{cut}$, by Lemma 17, $\mathbb{P}(\mathsf{tail}_3(\mathcal{W}_{k,n}) \leq m_{cut}) \to 0$. Furthermore, $K_m$ is bounded since $\lim_{m \to +\infty} K_m = 1$ according to Lemma 18. It follows that:

$$\sum_{m \leq m_{cut}} p_2 K_m \mathbb{P}(\mathsf{tail}_3(\mathcal{W}_{k,n}) = m) \to 0.$$

- Take $n$ sufficiently large that $m_{cut} \geq M$ (where $M$ is defined in Lemma 18). Then, for any $m > m_{cut}$, $1 \leq K_m \leq 1 + \frac{\gamma}{m_{cut}}$.

Hence, we obtain:

$$\mathbb{P}(\max \mathcal{W}_{k,n} \bmod 3 = 0) - p_2 \leq \sum_{m \leq m_{cut}} p_2 K_m \mathbb{P}(\mathsf{tail}_3(\mathcal{W}_{k,n}) = m) + p_2 \frac{\gamma}{(\beta n)^\alpha},$$

and the right-hand side tends to 0 as $n$ tends to infinity. By Proposition 9, we have proved that $\mu(C_3^t(x)_0 = 0) \to p_2$. □

Proof. (Theorem 6 for $a \neq 0$). The proof follows the case $a = 0$, the single difference being that the embedded symmetric walk corresponds to visits of $3\mathbb{Z} + a$ instead of $3\mathbb{Z} = 3\mathbb{Z} + 0$. Applying to the whole proof the map $(p_0, p_1, p_2) \longrightarrow (p_{(a+0) \bmod 3}, p_{(a+1) \bmod 3}, p_{(a+2) \bmod 3})$, we obtain that the asymptotic probability of the state $a$ is $p_{(a-1) \bmod 3}$ as expected. □

## 4. Conclusion.

*Other models.* It is natural to ask whether a similar phenomenon occurs in more complex models of cyclic dominance.

The easiest extension is to consider that each predator has a fixed probability $p < 1$ to replace its prey (probabilistic version). Although the global behaviour seems similar (see Figure 8), the probability that a small region surrounded by a predator and a prey disappears (instead of surviving with constant size) may impact the early dynamics; noise has been shown to create the possibility of this kind of "unlucky extinctions" in non-spatial models [26]. Experimental evidence seems to indicate a much slower convergence and no obvious numerical relationship between the initial parameters and asymptotic probability.

In higher dimension, we are so far from a complete understanding of the dynamics and limit measure [9] that no conjecture seems possible.

We believe that the random walk approach (Proposition 9) can be adapted to more general prey/predator relationship. Consider the prey/predator graph, where the oriented edge $i \to j$ means that $i$ is a predator for $j$. Beyond the simple 3-state cyclic dominance, more complex predator/prey graphs have been observed in nature [33, 31, 25]. As a first example, we believe Proposition 9 holds on alphabets of size $2k + 1$ where each state $n$ has $k$ predators $n + 1, n + 2, \ldots, n + k$ and $k$ preys $n - 1, \ldots, n - k$ (modulo $2k + 1$). The random walk has steps in $\{-k, +k\}$, with the same condition that $W[x]_i$ must be equal to $x_i$ modulo $2k + 1$.

The clearest limit to our approach is the presence of neutral particles, which can interact with other particles in ways that seem difficult to describe in terms of a simple height function (as an example, in the 4-state cyclic cellular automaton, a neutral particle can turn a positive particle into a negative particle or vice versa). The absence of neutral particle means that the prey/predator relationship corresponds to an orientation of the complete graph (a *tournament*).

*Max path preservation.* We believe that a necessary and sufficient condition for Proposition 9 to hold is the following. Assume the alphabet is $\mathbb{Z}/n\mathbb{Z}$ and denote $a < b < c$ if, by incrementing $a$ by one repeatedly, one reaches $b$ before $c$. A complete graph orientation is *max path preserving* if, for any triplet of distinct vertices/species $a < b < c$, if $a$ and $b$ are predators for $c$, ($a \to c$ and $b \to c$) then $b$ is a predator for $a$ ($b \to a$).

The definition of the walk $W[x]$ becomes the following:

- $w_0 = x_0$;
- if $x_i = x_{i+1}$ then $w_{i+1} = w_i$;
- if $x_i$ is a prey for $x_{i+1}$ then $w_{i+1}$ is the value equal to $x_{i+1}$ modulo $n$ in $\{x_i + 1, \ldots, x_i + n - 1\}$;
- if $x_i$ is a predator for $x_{i+1}$ then $w_{i+1}$ is the value equal to $x_{i+1}$ modulo $n$ in $\{x_i - 1, \ldots, x_i - (n - 1)\}$.

To understand the max path preserving assumption, consider the following situation: a factor $x_i x_{i+1} x_{i+2} = acb$ such that $a < b < c$ and $a$ and $b$ are predators for $c$. We have $W[x]_i > W[x]_{i+1} < W[x]_{i+2}$ and, since $a < b < c$, $W[x]_i < W[x]_{i+2}$. Since this factor becomes $ab$ at the next time step and the walk steps up, we want $b$ to predate $a$.

Brute force enumeration up to $n = 6$ suggests three families of prey-predator graphs for $n$ species with the max-path preserving property:

- the $n$ total orders compatible with the cyclic increments: $k, k-1, \ldots, 0, n-1, n-2, \ldots, k+1$ for $k$ from 0 to $n-1$ (the corresponding cellular automata is uninteresting as $k-1$ dominates every other state).
- some strongly connected prey/predator graphs where $0, n-1, n-2, \ldots, 1, 0$ forms a Hamiltonian cycle.
- some strongly connected prey/predator graphs where $0, n-1, n-2, \ldots, 1, 0$ is not a Hamiltonian cycle (but there is at least one Hamiltonian cycle, like for any strongly connected tournament).

The smallest example of the third family corresponds to the 3-state cyclic automaton where prey/predator relations are reversed. Notice that the resulting walk on $\mathbb{Z}$ consists of steps $\pm 2$ instead of steps $\pm 1$. Of course, in this case, we could consider the minimum on the walk on steps $\pm 1$ instead; but in general, we do not know if relabelling the prey/predator graphs of the third family according to one of the possibly many Hamiltonian cycles leads to a structure similar to the second family.

Up to $n = 6$, the last two families are counted by the Eulerian numbers (**A000295** in OEIS):

$$e_n := \sum_{k=1}^{\lfloor (n-1)/2 \rfloor} \binom{n}{2k+1}.$$

We conjecture the following characterisation: any orientation in the second family is defined by selecting an odd number of vertices $0 \leq c_1 < c_2 < \ldots < c_{2k+1} \leq n-1$, where $3 \leq 2k+1 \leq n$ and the order corresponds to the numbering given by the Hamiltonian cycle $0, n-1, n-2, \ldots, 1, 0$. The convention is that for any state $j$ in the cyclic interval $]c_i, c_{i+1}]$, $j$ is a predator for any state in $[c_{i-k \bmod 2k+1}, j[$. Since the number of possible cycles of $2k+1$ vertices is counted by the binomial coefficient $\binom{n}{2k+1}$, this characterisation would imply the counting above. Empirically, all orientations defined by this conjectural characterisation satisfy the max path preservation up to $n = 11$ species.

If Proposition 9 does generalise to these three graph families, it remains to check that the rest of the proof does as well. We expect a more systematic analysis using functional equations.

*Numerical results.* In order to conjecture a relationship between the asymptotic probability of each state and the parameters of the initial Bernoulli measure, we performed numerical simulations of various cellular automata: the 3-state cyclic; the 3-state cyclic with probability (invasion rate) 1/2; the 4-state cyclic; and the two cellular automata corresponding to the predator/prey graphs represented below.

$G_1$ :

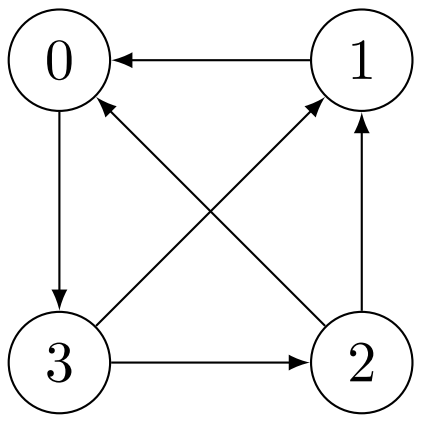

$G_2$ :

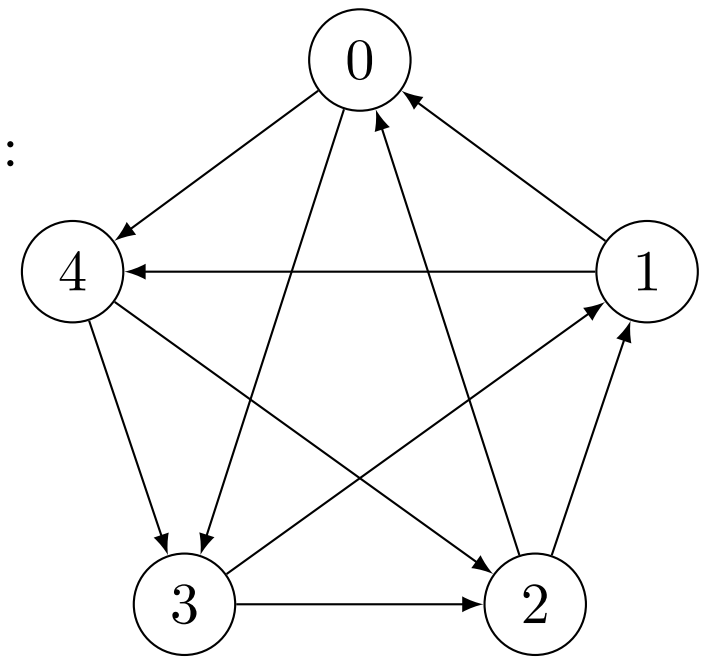

These results do not suggest any clear conjectural relationship, but we include them for possible future work.

For each cellular automaton $F$, we fixed values for the parameters of the initial Bernoulli measure that are distinct enough to be clearly distinguished.

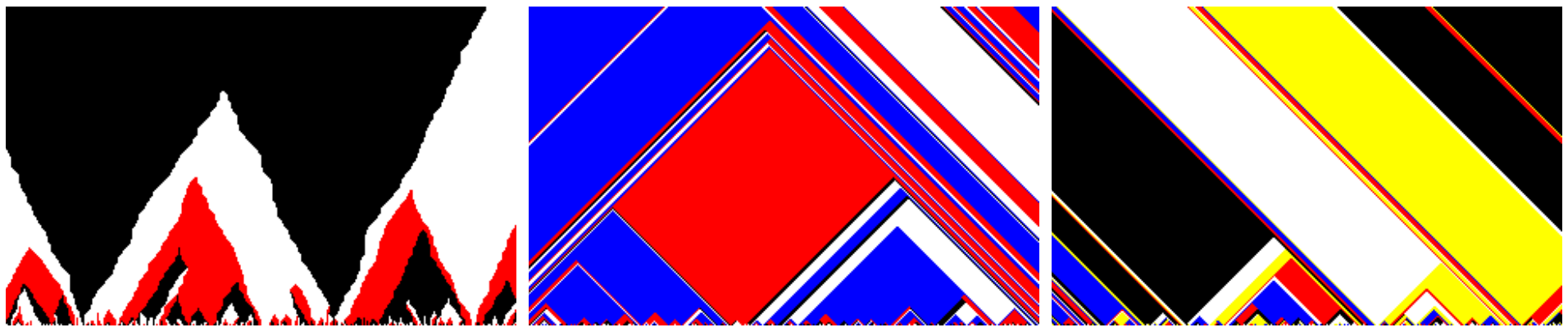

*Fig 8. Left to right, the 3-state cyclic cellular automaton with invasion rate* 1/2*, and the cellular automata corresponding to graphs* $G_1$ *and* $G_2$*.*

We computed the values of $(F^t(x)_0)_{0\leq t\leq 150}$ for 100000 random configurations $x$. The error margin (95% confidence) is $\pm 0.002$ for all the values.

| | | 3-state cyclic CA | | | 3-state w/ proba. 1/2 | | |
|---|---|---|---|---|---|---|---|
| states | | 0 | 1 | 2 | 0 | 1 | 2 |
| parameters | | 0.1 | 0.3 | 0.6 | 0.1 | 0.3 | 0.6 |
| time | 0 | 0.100 | 0.301 | 0.599 | 0.100 | 0.300 | 0.600 |
| | 50 | 0.550 | 0.087 | 0.362 | 0.760 | 0.087 | 0.153 |
| | 100 | 0.565 | 0.090 | 0.345 | 0.780 | 0.147 | 0.074 |
| | 150 | 0.572 | 0.091 | 0.336 | 0.742 | 0.198 | 0.060 |

| | | 4-state cyclic CA | | | | 4-state CA (graph $G_1$) | | | |
|---|---|---|---|---|---|---|---|---|---|
| states | | 0 | 1 | 2 | 3 | 0 | 1 | 2 | 3 |
| parameters | | 0.05 | 0.15 | 0.3 | 0.5 | 0.05 | 0.15 | 0.3 | 0.5 |
| time | 0 | 0.050 | 0.149 | 0.300 | 0.500 | 0.051 | 0.149 | 0.302 | 0.498 |
| | 50 | 0.065 | 0.578 | 0.018 | 0.339 | 0.369 | 0.009 | 0.040 | 0.580 |
| | 100 | 0.122 | 0.615 | 0.019 | 0.354 | 0.407 | 0.009 | 0.040 | 0.544 |
| | 150 | 0.005 | 0.615 | 0.123 | 0.368 | 0.424 | 0.010 | 0.041 | 0.525 |

| | | 5-state CA (graph $G_2$) | | | | |
|---|---|---|---|---|---|---|
| states | | 0 | 1 | 2 | 3 | 4 |
| parameters | | 0.025 | 0.075 | 0.15 | 0.3 | 0.45 |
| time | 0 | 0.025 | 0.076 | 0.152 | 0.298 | 0.450 |
| | 50 | 0.169 | 0.354 | 0.026 | 0.061 | 0.390 |
| | 100 | 0.165 | 0.376 | 0.028 | 0.061 | 0.370 |
| | 150 | 0.165 | 0.385 | 0.027 | 0.061 | 0.362 |

**Acknowledgements.** This work benefits from discussions between Jean-François Marckert and the two authors. The first author is grateful towards Antony Quas for his help and remarks on an earlier version of the paper.

E-mail: hellouin@lri.fr

E-mail: yvan.leborgne@labri.fr